\input amstex
\UseAMSsymbols
\NoBlackBoxes
\documentstyle{amsppt}
\magnification=\magstep1
\vsize=7.5in
\topmatter
\title{The prime number race and zeros
 of $L$-functions off the critical line}
\endtitle
\author Kevin Ford$^*$, Sergei Konyagin$^\dag$
\endauthor
\thanks
$^*$ First author supported in part by National Science Foundation grant
DMS-0070618.
\endthanks
\thanks
$^\dag$ Second author supported by INTAS grant N.~99-01080.
\endthanks
\rightheadtext{Prime number race and zeros of $L$-functions}
\leftheadtext{K. Ford, S. Konyagin}
\date April 27, 2002 \enddate
\address\noindent K.F.: Department of Mathematics,
 University of South Carolina, Columbia, SC 29208; \linebreak
S.K.: Department of Mechanics and Mathematics, Moscow State University,
Moscow 119899, Russia.
\endaddress
\abstract
We examine the effects of certain hypothetical
 configurations of zeros of Dirichlet $L$-functions lying off the critical
 line  on the distribution of primes in arithmetic progressions.
\endabstract
\endtopmatter
\document

\predefine\barunder{\b}
\def\ord{\text{ord}}
\def\flr#1{\left\lfloor #1 \right\rfloor}
\def\ceil#1{\left\lceil #1 \right\rceil}
\def\near#1{\left\| #1 \right\|}

\def\sg{\sigma}
\def\g{\gamma}
\def\del{\delta}
\redefine\b{\beta}
\define\lb{\lambda}
\def\lam{\lambda}
\define\Chi{\chi}
\define\bc{\overline{\chi}}

\def\curly{\Cal}
\def\BB{\curly B}

\redefine\le{\leqslant}
\redefine\ge{\geqslant}

\def\a{\alpha}

\define\({\left(}
\define\){\right)}
\define\pfrac#1#2{\( \frac{#1}{#2} \)}
\define\bfrac#1#2{\left[ \frac{#1}{#2} \right]}

\head 1. Introduction \endhead

Let $\pi_{q,a}(x)$ denote the number of primes $p\le x$ with $p\equiv
a\pmod{q}$.  The study of the relative magnitudes of the functions
$\pi_{q,a}(x)$ for a fixed $q$ and varying $a$ is known colloquially as
the ``prime race problem'' or ``Shanks-R\'enyi prime race problem''.
Fix $q$ and distinct residues $a_1,\ldots,a_r$ with $(a_i,q)=1$ for each $i$.
As colorfully described in the first paper of [KT1], consider a game
with $r$ players
called ``$1$'' through ``$r$'', and at time $t$, each player ``$j$'' has
a score of $\pi_{q,a_j}(t)$ (i.e. player ``$j$'' scores 1 point whenever
$t$ reaches a prime $\equiv a_j \pmod{q}$).  As $t\to \infty$, will
each player take the lead infinitely often?  More generally, 
will all $r!$ orderings of the players occur for infinitely many integers $t$?
It is generally believed that the answers to both questions is yes, for
all $q,a_1,\ldots,a_r$.

As first noted by Chebyshev [Ch] in 1853, some orderings
may occur far less frequently than others (e.g. if $q=3$, $a_1=1$, $a_2=2$,
then player ``$1$'' takes the lead for the first time when 
$t=608,981,813,029$ [BH]).  More generally, when $r=2$, $a_1$ is a 
quadratic residue modulo $q$, and $a_2$ is a quadratic non-residue modulo $q$,
$\pi_{q,a_2}(x)-\pi_{q,a_1}(x)$ tends to be positive more often
than it is negative (this phenomenon is now called ``Chebyshev's bias'').
In 1914, Littlewood [L] proved that both functions
$\pi_{4,3}(x)-\pi_{4,1}(x)$ and
$\pi_{3,2}(x)-\pi_{3,1}(x)$ change sign infinitely often.
Later Knapowski and Tur\`an ([KT1], [KT2])
 proved for many $q,a,b$ that $\pi_{q,b}(x)-\pi_{q,a}(x)$
changes sign infinitely often.
The distribution of the functions $\pi_{q,a}(x)$ is closely
linked with the distribution of the zeros in the critical
strip $0<\Re s < 1$ of the Dirichlet $L$-functions $L(s,\chi)$
for the characters $\chi$ modulo $q$. 
Some of the results of Knapowski and Tur\`an are proved under the
assumption that the functions $L(s,\chi)$ have no real zeros in $(0,1)$,
or that for some number $K_q$, the zeros of the functions $L(s,\chi)$  
with $|\Im s| \le K_q$ all have real part equal to $\frac12$.

Theoretical results for $r>2$ are more scant, all depending on the
unproven Extended Riemann Hypothesis for 
$q$ (abbreviated ERH$_q$), which
states that all these zeros lie on the critical line $\Re s = \frac12$. 
  Kaczorowski ([K1], [K2], [K3])
has shown that the truth of ERH$_q$ implies that for many
$r$-tuples $(q,a_1,\cdots,a_r)$,  $\pi_{q,a_1}(x) > \cdots
 > \pi_{q,a_r}(x)$ for arbitrarily large $x$.
If, in addition to ERH$_q$, one assumes that
the collection of non-trivial zeros of the $L$-functions for 
characters modulo $q$ are linearly independent over the rationals (GSH$_q$,
the grand simplicity hypothesis) ,
 Rubinstein and Sarnak 
[RS] have shown that for any $r$-tuple of coprime residue
classes $a_1,\ldots,a_r$ modulo $q$, that all $r!$ orderings of the
functions $\pi_{q,a_i}(x)$ occur for infinitely many integers $x$.  In fact
they prove more, that the logarithmic density of the set of real $x$ for which 
any such inequality occurs exists and is positive.

In light of the results of Littlewood  and of Knapowski and
Tur\`an, one may ask if such results for $r>2$ may be proved without the
assumption of ERH$_q$.
In particular, can it be shown, for some quadruples $(q,a_1,a_2,a_3)$,
that the 6 orderings of the functions $\pi_{q,a_i}(x)$
occur for infinitely many integers $x$,
without the assumption of ERH$_q$ (while still
allowing the assumption that zeros with imaginary part $< K_q$
lie on the critical line for some constant $K_q$)?  In this paper we
answer this question in the negative (in a sense) for all quadruples
$(q,a_1,a_2,a_3)$.  Thus, in a sense the hypothesis
ERH$_q$ is a necessary condition for proving
any such results when $r>2$.

Let $C_q$ be the set of non-principal characters modulo $q$.
Let $D=(q,a_1,a_2,a_3)$, where $a_1,a_2,a_3$
 are distinct residues
modulo $q$ which are coprime to $q$.
Suppose for each $\chi\in C_q$, $B(\chi)$ is a sequence of complex numbers
with positive imaginary part (possibly empty, duplicates allowed), and
denote by $\BB$ the system of $B(\chi)$ for $\chi\in C_q$.  Let $n(\rho,\chi)$
be the number of occurrences of the number $\rho$ in $B(\chi)$.
The system $\BB$ is called a {\it barrier} for $D$ if the following hold:

\item{(i)} all numbers in each $B(\chi)$ have real part in $[\beta_2,\beta_3]$,
where $\frac12 < \beta_2 < \beta_3\le 1$; 

\item{(ii)} for some $\beta_1$ satisfying $\frac12 \le \beta_1 < \beta_2$,
if we assume that
for each $\chi\in C_q$ and $\rho\in B(\chi)$, $L(s,\chi)$ has a zero
of multiplicity $n(\rho,\chi)$ at $s=\rho$, and all other zeros of $L(s,\chi)$
in the upper half plane have real part $\le \beta_1$, then one of the six 
orderings of the 
three functions $\pi_{q,a_i}(x)$ does not occur for large $x$.

If each sequence $B(\chi)$ is finite, we call $\BB$ a {\it finite
barrier} for $D$ and denote by $|\BB|$ the sum of the number of elements of each
sequence $B(\chi)$, counted according to multiplicity.

\proclaim{Theorem 1}
For every real numbers $\tau>0$ and $\sigma>\frac12$ and
every $D=(q,a_1,a_2,a_3)$, there is a finite
barrier for $D$, where each sequence
$B(\chi)$ consists of numbers with real part $\le \sigma$ and imaginary
part $>\tau$.  In fact, for most $D$, there is a barrier with $|\BB|\le 3$.
\endproclaim

We do not claim that the falsity of ERH$_q$ implies that one of the six
orderings does not occur for large $x$.  For example, take $\sg > \frac12$,
and suppose each non-principal character modulo  $q$ has a unique zero
with positive imaginary part to the right of the critical line, at
$\sg + i \gamma_\chi$.  If the numbers $\gamma_\chi$ are linearly independent
over the rationals, it follows easily from Lemma 1.1 below and the
Kronecker-Weyl Theorem that in fact
all $\phi(q)!$ orderings of the functions $\{ \pi_{q,a}(x) : (a,q)=1 \}$
occur for an unbounded set of $x$.

We now present a general formula for $\pi_{q,a}(x)$ in terms of the zeros
of the functions $L(s,\chi)$. 
Throughout this paper, constants implied by the Landau $O-$ and 
Vinogradov $\ll-$ symbols may depend on $q$, but not on any
other variable.

\proclaim{Lemma 1.1} Let $\beta\ge \frac12$, $x\ge 10$
and for each $\chi\in C_q$,
let $B(\chi)$ be the sequence of zeros (duplicates allowed)
of $L(s,\chi)$ with  $\Re s>\beta$ and
$\Im s > 0$.   Suppose further that
all $L(s,\chi)$ are zero-free on the real segment $0<s<1$.  If
$(a,q)=(b,q)=1$ and $x$ is sufficiently large, then
$$
\phi(q)\( \pi_{q,a}(x) - \pi_{q,b}(x) \) = -2\Re \left[ \sum_{\chi\in C_q} 
(\bc(a)-\bc(b)) \sum_{\Sb \rho\in B(\chi) \\ |\Im \rho|\le x \endSb}
 f(\rho) \right] + O(x^\b \log^2 x),
$$
where
$$
f(\rho) := \frac{x^{\rho}}{\rho\log x} + \frac{1}{\rho}
\int_2^x \frac{t^{\rho}}{t\log^2 t}
\, dt= \frac{x^{\rho}}{\rho\log x} + O\pfrac{x^{\Re \rho}}{|\rho|^2\log^2 x}.
$$
\endproclaim

\demo{Proof}  Let $\Lambda(n)$ be the von Mangolt function, and define
$$
\Psi_{q,a}(x) = \sum_{\Sb n\le x \\ n\equiv a\pmod{q} \endSb} \Lambda(n),
\qquad \Psi(x;\Chi) = \sum_{n\le x} \Lambda(n) \chi(n).
$$
Let $D_q$ be the set of all Dirichlet characters $\chi$ modulo $q$.
Then
$$
\split
\pi_{q,a}(x) &= \sum_{\Sb n\le x \\ n\equiv a\pmod{q} \endSb}
 \frac{\Lambda(n)}{\log n} + O(x^{1/2}) \\
&= \int_{2^-}^x \frac{d\Psi_{q,a}(t)}{\log t} + O(x^{1/2}) \\
&= \frac{\Psi_{q,a}(x)}{\log x} + \int_2^x \frac{\Psi_{q,a}(t)}{t\log^2 t}\,
    dt + O(x^{1/2})\\
&= \frac{1}{\phi(q)}
   \sum_{\chi\in D_q} \bc(a) \( \frac{\Psi(x;\chi)}{\log x} + \int_2^x
   \frac{\Psi(t;\chi)}{t\log^2 t}\, dt \) + O(x^{1/2}).
\endsplit
$$
Then
$$
\phi(q)(\pi_{q,a}(x)-\pi_{q,b}(x)) = \sum_{\chi\in C_q} (\bc(a)-\bc(b)) 
   \( \frac{\Psi(x;\chi)}{\log x} + \int_2^x
   \frac{\Psi(t;\chi)}{t\log^2 t}\, dt \) + O(x^{1/2}). \tag{1.1}
$$
By well-know explicit formulas (Ch. 19, (7)--(8) in [D]),
when $\chi\in C_q$,
$$
\Psi(x;\chi) = - \sum_{|\Im \rho| \le x} \frac{x^\rho}{\rho} + O\(\log^2 x \),
\tag{1.2}
$$
where the sum is over zeros $\rho$ of $L(s,\chi)$ with $0<\Re \rho < 1$.
Since the number of zeros
with $0\le \Im \rho \le T$ is $O(T\log T)$ ([D], Ch. 16, (1)), by partial
summation we have
$$
\sum_{\Sb 0< \Im \rho \le x \\ \Re \rho \le \beta \endSb}
  \left| \frac{x^{\rho}}{\rho} \right| \le
  x^\b \sum_{0<\Im \rho \le x} \frac{1}{|\rho|} \ll x^{\beta}\log^2 x.
$$
The implied constant depends on the character, and hence only on $q$.
By (1.2),
$$
\Psi(x;\chi) = - \sum_{\Sb |\Im \rho| \le x \\ \Re \rho > \b \endSb}
\frac{x^\rho}{\rho} + O\(x^{\beta}\log^2 x \).
\tag{1.3}
$$
The first part of the lemma follows by inserting (1.3) into (1.1) and
combining zeros $\rho$ of $L(s,\chi)$ and
$\overline{\rho}$ of $L(s,\bc)$.
Lastly, if $\frac12 \le \sg =\Re \rho$, integration by parts gives
$$
\split
\int_2^x \frac{t^\rho}{t\log^2 t}\, dt &= \left. \frac{t^\rho}{\rho \log^2 t}
  \right|_2^x + \frac{2}{\rho} \int_2^x \frac{t^{\rho-1}}{\log^3 t}\, dt \\
&\ll \frac{x^{\sigma}}{|\rho|\log^2 x} + \frac{1}{|\rho|} \left[\frac{1}
  {\log^3 2} \int_2^{\sqrt{x}} t^{\sigma-1}\, dt + \frac{8}{\log^3 x} 
  \int_{\sqrt{x}}^{x} t^{\sigma-1}\, dt \right] \\
&\ll \frac{x^{\sigma}}{|\rho|\log^2 x}.
\endsplit
$$
This completes the proof of the lemma.
\qed\enddemo

In the next three sections, we show several methods for constructing 
barriers, which, by Lemma 1.1, boils down to analyzing the two functions
$$
\Re \sum_{\chi\in C_q} (\bc(a_j)-\bc(a_3)) \sum_{\rho\in B(\chi)} 
\frac{x^{\rho}}{\rho} \quad (j=1,2).
$$
In section 2 we construct a barrier using two simple zeros (one of which
may be a zero for several characters).  Section 3 details a method using a
zero for $L(s,\chi)$ and a zero for $L(s,\chi^2)$ (for most $D$ these
are simple or double zeros). 
Lastly, section 4 presents a more general method with two numbers,
which are zeros for each character of certain high multiplicities.
Together, the three constructions provide barriers for all 
quadruples $(q,a_1,a_2,a_3)$.

%  We also make sure that our constructions
%do no violate known zero-density bounds for $L$-functions ([Se], [XX]);
%some of the constructions assume that we have a number $\rho$ 
%which is a zero for many $\chi \in C_q$, possibly with
%high multiplicities.  To satisfy the density hypotheses, our 
%constructions of barriers $\BB$ will always satisfy
%$$
%|\BB| \ll q^{1/2}.
%$$

All of the constructions in sections 2--4 involve two
zeros, one with imaginary part $t$ and the other with imaginary part
$2t$.  Thus, we assume that both ERH$_q$ and GSH$_q$ are false. 
Answering a question posed by Peter Sarnak, in
section 5 we construct a barrier (with an infinite set $B(\chi)$)
where the imaginary parts of the numbers in the sets $B(\chi)$  are
linearly independent; in particular, we assume all zeros of each
$L(s,\chi)$ are simple, and $L(s,\chi_1)=0=L(s,\chi_2)$ does not
occur for $\chi_1 \ne \chi_2$ and $\Re s >\beta_2$.

We adopt the notations
 $e(z)=e^{2\pi i z}$, $\lfloor x \rfloor$ is the greatest integer
$\le x$, $\lceil x \rceil$ is the least integer $\ge x$,
 $\{ x\} = x - \flr{x}$ is the fractional
part of $x$, and $\near{x}$ is the distance from $x$ to the nearest integer.
Also,
$\arg z$ is the argument of the nonzero complex number $z$ lying in
$[-\pi,\pi)$.  Throughout, $q=5$ or $q\ge 7$, and $(a_1,q)=(a_2,q)=
(a_3,q)=1$.

%%%%%%%%%%%%%%%%%%%%%%%%%%%%%%%%%%%%%%%%%%%%%%%%%%%%%%%%%%%%%%%%%%%%%%%%%%%
%
\head 2.  First construction \endhead
%
%%%%%%%%%%%%%%%%%%%%%%%%%%%%%%%%%%%%%%%%%%%%%%%%%%%%%%%%%%%%%%%%%%%%%%%%%%%
%
%%%%%%%%%%%%%%%%%%%%%%%%%%%%%%%%%%%%%%%%%%%%%%%%%%%%%%%%%%%%%%%%%%
%
%  Example from April 1999
%
%  (This is a generalization of Examples 1, 2 and 3, and
%  Corollary 1 from S.K. notes, April 1999). 
%
\proclaim{Lemma 2.1}
If, for some relabelling of the numbers
$a_i$, there is a set $S$ of nonprincipal
 Dirichlet characters modulo $q$  such that
$$
\sum_{\chi \in S} \chi(a_1) = \sum_{\chi \in S} \chi(a_2) \ne 
\sum_{\chi \in S} \chi(a_3),
$$
then there is a barrier $\BB$ for $D=(q,a_1,a_2,a_3)$ with $|\BB| \le |S|+1$.
\endproclaim

{\bf Remark.}  The hypotheses of Lemma 2.1 are satisfied when, for 
example, $q$ has a primitive root $g$, and $a_3/a_2$ is not in the subgroup
of $(\Bbb Z/q\Bbb Z)^*$ generated by $a_2/a_1$.  Writing $a_2/a_1\equiv g^f$,
we take the character with $\chi(g)=e(1/(f,\phi(q)))$ and $S=\{ \chi \}$.
\bigskip

\demo{Proof} Suppose $1/2 \le \b < \sigma_2 < \sigma_1 \le \min(\sigma,
0.501)$,
and let $\chi_2$ be a character with $\chi_2(a_1) \ne \chi_2(a_2)$
($\chi_2$ may or may not be in $S$).  Let
$T_q$ be a large number, depending only on $q$.
Let $\rho_1 = \sigma_1+it$, $\rho_2 = \sigma_2+2it$ where
$t > T_q$.  Suppose
$L(s,\chi)$ has a simple zero at $s=\rho_1$ for each $\chi\in S$, 
$L(s,\chi_2)$ has a simple zero at $s=\rho_2$, and no other non-trivial
zeros of any $L$-function in $C_q$ have real part exceeding $\b$.
Let 
$$
D_1(x):=\phi(q)(\pi_{q,a_1}(x) - \pi_{q,a_2}(x)), \qquad
D_2(x):=\phi(q)(\pi_{q,a_3}(x) - \pi_{q,a_2}(x)).
$$
By Lemma 1.1 and our hypotheses, if $x$ is sufficiently large,
$$
\split
D_1(x) &= \frac{2 x^{\sigma_2}}{\log x} \left[
\Re \(\frac{e^{2it\log x}}{\sigma_2 + 2it} W\) + O\pfrac{1}{\log x} \right],
\quad W=\bc_2(a_2)-\bc_2(a_1), \\
D_2(x) &= \frac{2 x^{\sigma_1}}{\log x} \left[
\Re \( \frac{e^{it\log x}}{\sigma_1 + it} Z\) + O\pfrac{1}{\log x} \right],
\quad Z=\sum_{\chi\in S} \( \bc(a_2)-\bc(a_3) \).
\endsplit
$$
Define
$$
\split
A(x) &= \near{\frac{1}{\pi} \arg \( \frac{e^{it\log x}}{\sigma_1 + it} Z\) - 
\frac12} \\
&=\near{\frac{1}{\pi}(t\log x+\arg Z + \tan^{-1} (\sigma_1/t))}.
\endsplit
$$
If $A(x) \ge (\log x)^{-1/2}$, then
$|D_2(x)| \gg x^{\sigma_1}/\log^{3/2} x$.
But $D_1(x)=O(x^{\sigma_2})$, so for such $x$,
 $\pi_{q,a_3}(x)$ is either the largest
or the smallest of the three functions.  When $A(x) < (\log x)^{-1/2}$, then
$$
\split
C(x) &:= \arg \( \frac{e^{2it\log x}}{\sigma_2 + 2it} W\) \\
&\equiv \arg W - \frac{\pi}2 + \tan^{-1}\pfrac{\sigma_2}{2t} +2t\log x \\
&\equiv \arg W + \tan^{-1}\pfrac{\sigma_2}{2t}- 2\arg Z -
   2 \tan^{-1} \pfrac{\sigma_1}{t} + O\pfrac{1}{\sqrt{\log x}} \\
&\equiv  \arg W - 2 \arg Z - F(x) \pmod{\pi},
\endsplit
$$
where $1/(2t) < F(x) < 1/t$ for large $x$.
The number of possibilities for
$\arg W - 2 \arg Z$ depends only on $q$, hence we may assume either 
$$
B=\left\{ \frac{1}{\pi}(\arg W-2\arg Z)\right\} - \frac12
$$
 satisfies either $B=0$ or $|B|>2/t \ge 2 F(x)$ (by taking $T_q$ sufficiently
large).  We have
$$
C(x) \equiv \pi B + \frac{\pi}2 - F(x) \pmod{\pi}.
$$
 If $B=0$, then $C(x)$
is either $\pi/2 - F(x)$ or $3\pi/2-F(x)$ (mod $2\pi$), whence
$D_1(x)$ takes only one sign for such $x$.
 Likewise, 
$C(x) \in (\pi/2+2/t,\pi)$ if $B>2/t$ and
$C(x) \in (-F(x),\pi/2-2/t)$ if $B<-2/t$.
In all cases, when $A(x) < (\log x)^{-1/2}$,  $D_1(x)$
takes only one sign.  Therefore,
one of the orderings $\pi_{q,a_1}(x) > \pi_{q,a_3}(x) > \pi_{q,a_2}(x)$ or
 $\pi_{q,a_2}(x) > \pi_{q,a_3}(x) > \pi_{q,a_1}(x)$ does not
occur for large $x$. 
\qed\enddemo

{\bf Remark.}  By similar reasoning, for any integer $k\ge 2$ 
one may construct a barrier with one zero having imaginary part $t$ and
another zero having imaginary part $kt$.
\medskip

%%%%%%%%%%%%%%%%%%%%%%%%%%%%%%%%%%%%%%%%%%%%%%%%%%%%%%%%%%%%%%%%%%%%%%
%
%
\head 3. Second construction \endhead
%
%
%%%%%%%%%%%%%%%%%%%%%%%%%%%%%%%%%%%%%%%%%%%%%%%%%%%%%%%%%%%%%%%%%%%%%%

The basic idea of this section is to find a character $\chi$ so that
the values $\chi(a_1)$, $\chi(a_2)$, $\chi(a_3)$ are nicely 
spaced around the unit circle, but not too well spaced (e.g. cube roots
of 1 or translates thereof).  In almost all circumstances we can find
such a character.

\proclaim{Lemma 3.1}  Let $s_1=\ord_q(a_2/a_1)$, $s_2=\ord_q(a_3/a_2)$ and
$s_3=\ord_q(a_1/a_3)$.  If one of $s_1,s_2,s_3$ is not in
$\{ 3, 7, 13, 21\}$, then for some 
relabeling of the $a_i$'s, there is a Dirichlet character $\chi$ satisfying 
either
\roster
\item"(i)" $\chi(a_1) = \chi(a_2) \ne \chi(a_3)$; or
\item"(ii)" $\chi(a_i)=e(r_i)$ with $0 \le r_1 < r_2 < r_3 < 1$, and
$d_1=r_2-r_1$, $d_2=r_3-r_2$ satisfy
$$
\frac13 < d_1 \le d_2 < \frac12, \quad \text{ or } \quad
(d_1,d_2) \in \left\{ \( \tfrac{6}{19},\tfrac{9}{19} \), \( \tfrac{12}{37},
\tfrac{16}{37}\) \right\}.
\tag{3.1}
$$
\endroster
\endproclaim

\noindent
{\bf Remark.}
In the case that (i) holds, the hypotheses of Lemma 2.1 hold with
$S=\{ \chi \}$, and thus there is a finite barrier for $D$ with
$|\BB|=2$.  Therefore, in this section we confine ourselves with the
case that (ii) holds (Lemma 3.5 below).

Before proving Lemma 3.1, we begin with some simple lemmas about the
existence of characters with certain properties.

\proclaim{Lemma 3.2}  Suppose $q\ge 3$ and $(b,q)=1$.  Let $m$ be the order of
$b$ modulo $q$. 
 Then there is a Dirichlet
character $\chi$ modulo $q$ with $\chi(b) = e(1/m)$.
\endproclaim

\demo{Proof} Suppose $g_1,\ldots,g_t$ generate $(\Bbb Z/q\Bbb Z)^*$ and
$b=g_1^{f_1} \cdots g_t^{f_t}$.  Let $s_i = \ord_q g_i$ for each $i$,
and $s_i'$ be the order of $g_i^{f_i}$.  Then $s_i' = s_i/(f_i,s_i)$
and $m=\text{lcm}[s_1',\ldots,s_t']$.
Let $f_i' = f_i/(f_i,s_i)$, so in particular $(s_i',f_i')=1$.  The
$\gcd$ of the $t+1$ numbers $m,f_i'm/s_i'$ is 1, so there are integers
$h_1,\ldots, h_t$ so that
$\sum h_i \frac{f_i'm}{s_i'} \equiv 1 \pmod{m}.$
Take the character $\chi$ with $\chi(g_i)=e(h_i/s_i)$ for each $i$, then
$\chi(b) = \prod \chi(g_i)^{f_i} = e(h_1f_1'/s_1' + \cdots + 
h_t f_t'/s_t') = e(1/m)$.
\qed\enddemo

\proclaim{Lemma 3.3} Suppose $b,c$ are distinct residues modulo $q$
with $(b,q)=(c,q)=1$.  Suppose that $r | \ord_q b$
and for every $p^a \| r$ with $a\ge 1$, $p^{a+1} \nmid \ord_q c$.
Then there is a Dirichlet character $\chi$ modulo $q$ such that
$$
\chi(b) = e(1/r), \qquad \chi(c)^r = 1.
$$
\endproclaim

\demo{Proof}
Let $s_1=\ord_q b$ and $s_2 = \ord_q c$.  By Lemma 3.2, there is a character
$\chi_1$ with $\chi_1(b)=e(1/s_1)$ and therefore a character $\chi_2$ with
$\chi_2(b)=e(1/r)$.  Since $c$ has order $s_2$, $\chi_2(c) = e(g/s_2)$
for some integer $g$.  Write $s_2 = v u$ where $(u,r)=1$ and $v|r$.
Define $x$ by $xu \equiv 1\pmod{r}$, and let
$\chi=\chi_2^{xu}$.  Then $\chi(b)=\chi_2(b)^{xu} = e(1/r)$ and
$\chi(c) = e(gxu/s_2) = e(gx/v) = e(gx(r/v)/r)$.
\qed\enddemo

\proclaim{Definition}  An odd number $m$ is ``good'' if for every $j$,
$1 \le j \le m-1$, there is a number $k$ such that
 among the points $(0,k/m,kj/m) \mod 1$, either two are equal (and not equal
to the third), or two of the three distances $d_1,d_2,d_3$ (with sum = 1)
 between the points satisfy (3.1).
\endproclaim

{\bf Remark.} To prove that a number $m$ is good, we need only to check
$2\le j\le (m+1)/2$, since for $j=1$
we take $k=1$, and if $k$ works for $j=j_0$ then the same $k$ works for
$j=m+1-j_0$.

\proclaim{Lemma 3.4}  Every odd prime $p$ except $p\in P=\{ 3, 7, 13\}$ is
good, and for $p\in P$, $p^2$ is good.  Also, the numbers 39, 91 and
273 are good.
\endproclaim

\demo{Proof}  A short computation implies that if $p\in P$, then
$p$ is not good, but $p^2$ is good.  Also, by a short computation, all other
odd primes $\le 83$ are good, as well as 39, 91 and 273.  The following $j$
values have no associated $k$-value: for $m=3$, $j=2$; for $m=7$, $j=3,5$;
 for $m=13$; $j=3,5,6,8,9,11$; for $m=21$, $j=5,17$.

Suppose that $m=p > 84$ is prime and write each product $kj=\ell p+r$ with
$0 \le r < p$.  We shall prove that for each $j \in [2,\frac{p+1}2]$,
there is a $k$ so that two of the three distances satisfy
$\frac13 < d_1 \le d_2 < \frac12$.
We now divide up the $j \in [2,\frac{p+1}2]$ into 9 cases:
\medskip

{\bf Case I.}  $j\in \{ 3, 5, 7, \frac{p+1}{2} \}$.
 For
$j=3$ take $p/6 < k < 2p/9$ and for $j=5,7$ take any $k$
with $p/(2j) < k < p/(2j-2)$.  There is such a $k$ when $p>84$.
Then $p-jk$ and $jk-k$ both lie in $(p/3,p/2)$.
For $j=\frac{p+1}2$ take $k=2\ceil{p/3}$, then $r=\ceil{p/3}$, so both
$r$ and $k-r$ lie in $(p/3,p/2)$ for $p>6$.
% For $j=\frac{p+4}3$, take
% $k=3\ceil{p/9}+1$, so $r=\frac{p+1}3+\ceil{p/9}$, and for $j=\frac{p+5}{3}$
% take $k=\ceil{p/9}$, so $r=5\ceil{p/9}$.
\medskip

{\bf Case II.} $9 \le j < p/6+1$.  Take $m= \flr{\frac5{12}(j-1)}$.
Then 
$$
\frac{m+1/2}{j-1} \le \frac5{12} + \frac{1/2}{j-1} < \frac12, \qquad
\frac{m+1/3}{j-1} \ge \frac{5}{12} - \frac{7/12}{j-1} > 1/3.
$$
Therefore, if
$$
\frac{p(m+1/3)}{j-1} < k < \frac{p(m+1/2)}{j-1},
$$
then $k$ and $r-k$ lie in $(p/3,p/2)$.
But the above interval has length $p/(6j-6) > 1$, so such a $k$
exists.
\medskip

{\bf Case III.}  $2\le j < p/3+1, j$ even.  Take $k=\frac{p-1}2$.  Then
$r=p-j/2$ and both $k$ and $r-k$ lie in $(p/3,p/2)$.
\medskip

{\bf Case IV.} $p/3+1 < j < 3p/7, j$ even.  Take $h$ so that $1\le h
< \frac{p-3}{18}$ and
$$
\frac{2h+2/3}{6h+1}p < j < \frac{2h+1}{6h+1}p.
$$
The largest admissible $h$ is at least $\frac{p-19}{18}$, so
the above intervals cover $(\frac{p(p-13)}{3(p-16)}, 3p/7)$, which contains
$[ \frac{p+4}{3}, 3p/7)$ for $p>64$.  Then take $k=\frac{p-1}2 - 3h$,
so that $r \in (p/2, 2p/3)$.

\medskip

{\bf Case V.} $2p/5 +1 < j \le \frac{p-1}2, j$ even.  We take $h$ so that
$1\le h < \frac{p-3}{12}$ and
$$
\frac{2h}{4h+1}p < j-1 < \frac{2h+1/3}{4h+1}p.
$$
The largest admissible $h$ is at least $\frac{p-13}{12}$, so these intervals
cover $(2p/5,\frac{p(p-11)}{2(p-10)})$, which includes $(2p/5,\frac{p-3}{2}]$
for $p>13$.  Then take $k=\frac{p-1}{2}-2h$, so $r-k \in (p/3,p/2)$.
\medskip

{\bf Case VI.} $p/3+1 < j \le \frac{p-1}{2}, j$ odd.  Take $h$, $0 \le h
< \frac{p-15}{12}$ so that
$$
\frac{2h+1}{4h+3}p < j-1 < \frac{2h+4/3}{4h+3}p.
$$
Than take $k=\frac{p-3}2 - 2h$, so that $k-r \in (p/3,p/2)$.
The above intervals cover $(p/3,\frac{p-3}2]$
provided that $p>24$.
\medskip

{\bf Case VII.} $p/3 -1 < j < p/3+1$.
Write $j=\frac{p+t}{3}$, where
$-2 \le t\le 2$, $t\ne 0$.  Here we take $k=3\ceil{p/9} + b$,
where $0\le b\le 2$ and $t+3b\equiv w \pmod{9}$, $w\in \{5, 7\}$.
If $p>28$ then $k\in (p/3,p/2)$.  If $w=5$, then $r=5p/9 + E$,
where $|E| \le 22/9$.  Thus, $r\in (p/2,2p/3)$ when 
$p>44$.  When $w=7$, $t=1,b=2$, then $r\in (7p/9,7p/9+14/9]$.

\medskip

{\bf Case VIII.} $5p/21 < j < p/3-1, j$ odd.  Take $1\le h < \frac{p-3}{18}$
so that
$$
\frac{6h-1}{18h+3}p < j < \frac{2h}{6h+1}p.
$$
Take $k=\frac{p-1}{2} - 3h$, so $r\in (p/3,p/2)$.
The above intervals cover $(5p/21,p/3-1)$.
\medskip

{\bf Case IX}.  $p/6+1 < j < 5p/21, j$ odd.  If $p/5< j-1 < 4p/15$, take
$k=\frac{p-5}{2}$, so that $r \in (5p/6-5/2, p-5/2)$.  If
$p/7 < j-1 < 4p/21$, then $k=\frac{p-7}2$ works and if $5p/27 < j < 2p/9$
then $k=\frac{p-9}2$ works.
\qed\enddemo

\bigskip
\demo{Proof of Lemma 3.1}
By hypothesis, there are two possibilities: 
\item{(i)} some $s_i$ (say $s_1$) is divisible by a prime power $p^w$
other than $3, 7$, or $13$;

\item{(ii)} Each $s_i$ divides $273$ and some $s_i$
(say $s_1$) equals 39, 91 or 273.

Say $s_1$ is
divisible by $p^w$, with $p^{w+1} \nmid s_2$ and $p^{w+1} \nmid s_3$.  
By Lemma 3.3, there is a character $\chi_1$ with $\chi_1(a_2/a_1)=e(1/p^w)$
and $\chi_1(a_3/a_2) = e(m/p^w)$ for some integer $m$.  If $p=2$, let
$\chi = \chi_1^{2^{w-1}}$, so that $\chi(a_2/a_1)=-1$ and
$$
1 = \chi(a_2/a_1) \chi(a_3/a_2) \chi(a_1/a_3) = - \chi(a_3/a_2) \chi(a_1/a_3).
$$
But each character value on the right is either -1 or 1, so either
$\chi(a_2)=\chi(a_3)$ or $\chi(a_1)=\chi(a_3)$ and (i) is satisfied.
If $p$ is odd, let 
$\chi_2=\chi_1^{p^{w-1}}$ if $p\not\in P$ and $\chi_2=\chi_1^{p^{w-2}}$ if
 $p\in P$.  Then $\chi_2(a_2/a_1) = e(1/p^u)$, where $u=2$ if $p\in P$
and $u=1$ otherwise.  Write $\chi_2(a_3/a_2) = e(j/p^u)$.  If $j=0$ then
$\chi_2(a_2)=\chi_2(a_3)$ and (i) is satisfied.  Otherwise,
since $p^u$ is good by Lemma 3.4, there is a number $k$ so that two of the
three distances of the points $(0,k/p^u,kj/p^u) \pmod{1}$ satisfy (3.1).
Taking $\chi=\chi_2^k$ gives (ii) for some relabeling
of the $a_i$'s.

In the case that each $s_i$ divides 273 and $s_1 \in \{ 39, 91, 273\}$,
by Lemma 3.3 there is a character $\chi_1$ with $\chi_1(a_2/a_1)=e(1/r)$
and $\chi_1(a_3/a_2)=e(g/r)$ for some integer $g$.  (here $r=s_1$).  Since
$r$ is good by Lemma 3.4, there is a $k$ such that two of the
three distances of the points $(0,k/r,kj/r) \pmod{1}$ satisfy (3.1).
Taking $\chi=\chi_1^k$ gives (ii) for some relabeling
of the $a_i$'s.
\qed\enddemo

%
%%%%%%%%%%%%%%%%%%%%%%%%%%%%%%%%%%%%%%%%%%%%%%%%%%%%%%%%%%%%%%
%
%
% Main result November 2-11, 1999
%
%%%%%%%%%%%%%%%%%%%%%%%%%%%%%%%%%%%%%%%%%%%%%%%%%%%%%%%%%%%%%%
%

\proclaim{Lemma 3.5}  Suppose that for some relabeling of $a_1,a_2,a_3$ and
some Dirichlet character $\chi$ modulo $q$,
$\chi(a_i)=e(r_i)$ with $0 \le r_1 < r_2 < r_3 \le 2$,
$d_1=r_2-r_1$ and $d_2=r_3-r_2$ and $(d_1,d_2)$ satisfies (3.1).
Then there is a finite barrier $\BB$ for $D=(q,a_1,a_2,a_3)$
with $|\BB| \le 14$.  If $d_1>\frac13$, then $|\BB| \le 3$.
\endproclaim

\demo{Proof}  For some $1/2 \le \beta <\a \le \sg$
 and large $\gamma$, suppose $L(s,\chi)$ has a zero at $s=\a+
i\gamma$ of order $c_1$, and $L(s,\chi^2)$ has a zero at $s=\a+2i\gamma$
of order $c_2$, where 
$$
(c_1,c_2) = \cases (1,2) & d_1 > \tfrac13 \\ (5,9) & d_1=\tfrac{6}{19} \\
(3,5) & d_1=\tfrac{12}{37}. \endcases
$$
Suppose all other non-trivial zeros of $L$-functions
modulo $q$ have real part $\le \beta$.  Let
$$
\split
D_1(x) &= \frac{\phi(q)\log x}{x^\a} (\pi_{q,a_2}(x)-\pi_{q,a_1}(x)), \\
D_2(x) &= \frac{\phi(q)\log x}{x^\a} (\pi_{q,a_3}(x)-\pi_{q,a_2}(x)).
\endsplit
$$
Let $u=\log x$.   For large $x$, Lemma 1.1 and the identity
$$
\sin(a-b)-\sin(a-c)=2\cos (a - \tfrac{b+c}{2}) \sin ( \tfrac{c-b}{2} )
$$
give
$$
\aligned
D_1(x) &=\frac{4}{\gamma} \sum_{\ell=1}^2 \frac{c_\ell}{\ell} 
  \sin(d_1 \ell \pi) \cos(\ell \gamma u -(r_1+r_2)\pi \ell)+ O(1/\gamma^2), \\
D_2(x) &=\frac{4}{\gamma} \sum_{\ell=1}^2 \frac{c_\ell}{\ell} \sin(d_2 \ell
  \pi) \cos(\ell \gamma u - (r_2+r_3)\pi \ell) + O(1/\gamma^2).
\endaligned\tag{3.2}
$$
For $j=1,2$ define
$$
\split
g_j(y) &= c_1 \sin(\pi d_j) \cos y + \frac{c_2}{2} \sin(2\pi d_j) \cos 2y\\
&= c_1 \sin(\pi d_j) \( \cos y + \tfrac{c_2}{c_1} \cos (\pi d_j) \cos 2y  \).
\endsplit
$$
Because $0 < d_j < 1/2$, $\cos \pi d_j$ and $\sin \pi d_j$ are both positive.
We claim that
$$
\min (g_1(\gamma u - (r_1+r_2)\pi), g_2(\gamma u - (r_2+r_3)\pi))
 < 0 \qquad (u\ge 0), \tag{3.3}
$$
which is equivalent to showing
$$
\min\bigl( g_1(y), g_2(y-\pi(d_1+d_2)) \bigr) < 0
$$
for all real $y$.
Since $g_1$ and $g_2$ are periodic and continuous, in fact the minimum
above is $\le -\delta$ for some $\delta>0$.
If $\gamma$ is large (depending on $\delta$), this implies that one of the two
functions on the left in (3.2) is negative for all large $x$.  Thus
for large $x$, $\pi_{q,a_3}(x) > \pi_{q,a_2}(x) > \pi_{q,a_1}(x)$ does
not occur.

\def\lb{\lambda}
To prove (3.3), we consider the one parameter family of functions
$h(y;\lb) = \cos y + \lb \cos(2y)$ for $0 < \lb < 1$.  These are all even 
functions, so it suffices to look at $0 \le y \le \pi$.  We have $h(y;\lb)$
positive for $0 \le y < v_\lb$ and negative for $v_\lb < y \le \pi$, where
$v_\lb = \cos^{-1} [\frac{1}{4\lb} (-1+\sqrt{8\lb^2+1})]$.
As a function of $\lb$, $v_\lb$ decreases from $\pi/2$ at $\lb=0$ to $\pi/3$ at
$\lb=1$.
For $i=1,2$, let $z_i= v_{\lb_i}$ for $\lb_i=(c_2/c_1) \cos \pi d_i$.
Since $\pi(d_1+d_2) < \pi$, (3.3) will follow from
$$
z_1 + z_2 < \pi(d_1+d_2). \tag{3.4}
$$
When $(d_1,d_2) \in \{ ( \tfrac{6}{19},\tfrac{9}{19} ), ( \tfrac{12}{37},
\tfrac{16}{37})\}$, (3.4) follows by direct calculation.  When
$\frac13 < d_1$, we have $c_1=1$, $c_2=2$ and
$\lb_j=2\cos \pi d_j$ ($j=1,2$).  We claim for $j=1,2$ that
$z_j <\pi d_j$, or equivalently $\cos z_j>\cos \pi d_j=\frac12 \lb_j$.
Since $0<\lb_j < 1$, 
$$
\cos z_j = \frac{\sqrt{8\lb_j^2+1}-1}{4\lb_j} >
\frac{\sqrt{4\lb_j^4+4\lb_j^2+1}-1}{4\lb_j} = \frac{\lb_j}{2},
$$
which proves (3.4) in this case as well.
\qed
\enddemo

Combining Lemmas 3.1 and 3.5 gives the following.

\proclaim{Corollary 3.6}  Let $s_1=\ord_q(a_2/a_1)$, $s_2=\ord_q(a_3/a_2)$ and
$s_3=\ord_q(a_1/a_3)$.  If one of $s_1,s_2,s_3$ is not in
$\{ 3, 7, 13, 21\}$, then there is a finite barrier $\BB$ for $D$
with $|\BB| \le 14$. 
\endproclaim
\bigskip

%%%%%%%%%%%%%%%%%%%%%%%%%%%%%%%%%%%%%%%%%%%%%%%%%%%%%%%%%%%%%%%%%%%
%
%
%
\head 4. Third construction \endhead
%
% from S.K. Feb 25, 2000; modified 3/11/2000; simplified 11/12/2000
%
%%%%%%%%%%%%%%%%%%%%%%%%%%%%%%%%%%%%%%%%%%%%%%%%%%%%%%%%%%%%%%%%%%%

Throughout this section, we assume that
$a_1,a_2,a_3$ do not satisfy the conditions of Lemma 2.1.

\proclaim{Lemma 4.1} Let
$\chi$ be a character modulo $q$ such that there are at least
two different values among $\chi(a_1)$, $\chi(a_2)$, $\chi(a_3)$.
Then the following hold:\newline
(a) $\chi(a_1)$, $\chi(a_2)$, $\chi(a_3)$ are distinct;\newline
(b) $\Re\chi(a_1)$, $\Re\chi(a_2)$, $\Re\chi(a_3)$ are distinct;\newline
(c) All the values $\chi(a_1)$, $\chi(a_2)$, $\chi(a_3)$ are not $\pm1$.
\newline
(d) $\chi$ has order $\ge 7$.
\endproclaim
\demo{Proof} (a) If this does not hold, the conditions of Lemma 2.1 hold
with $S=\{ \chi\}$.

(b) If $\chi(a_1)=\overline\chi(a_2)$, then, by (a),
$\Re\chi(a_3)\neq\Re\chi(a_1)$, and the conditions of Lemma 2.1 hold for
$S=\{\chi,\overline\chi\}$.

(c) If $\chi(a_3)=1$ and $k$ is the order of the character $\chi$,
then the conditions of Lemma 2.1 hold for
$S=\{\chi,\chi^2,\dots,\chi^{k-1}\}$. If $\chi(a_3)=-1$ and none of
$\chi(a_i)=1$, then
$\chi^2(a_3)=1\neq\chi^2(a_1)$, and the conditions of Lemma 2.1 hold
for $S=\{\chi^2,\chi^4,\dots,\chi^{2h-2}\}$ where $h$ is the order of $\chi^2$.

(d) This follows directly from (b) and (c).
\enddemo

\proclaim{Lemma 4.2}  
There exists a character $\chi$ modulo $q$ of order $\ge 7$
such that
$$
\Re (\chi(a_3)-\chi(a_2)) \Re (\chi^2(a_2)-\chi^2(a_1)) \ne
\Re (\chi(a_2)-\chi(a_1)) \Re (\chi^2(a_3)-\chi^2(a_2))
\tag{4.1}
$$
and for some integers $h,k$ with $1\le h< k\le 3$,
$$
\Im (\chi^h(a_3)-\chi^h(a_2)) \Im (\chi^k(a_2)-\chi^k(a_1)) \ne
\Im (\chi^h(a_2)-\chi^h(a_1)) \Im (\chi^k(a_3)-\chi^k(a_2)).
\tag{4.2}
$$
\endproclaim

\demo{Proof} Let $\chi$ be any character modulo $q$ such that
$\chi(a_2/a_1) \ne 1$.  By Lemma 4.1 (a),
the values $\chi(a_1)$, $\chi(a_2)$, $\chi(a_3)$ are distinct.
Denote $\chi(a_j)=e^{2\pi i \varphi_j}$ ($j=1,2,3$). By Lemma 4.1 (b),
the values $\cos(\varphi_1)$, $\cos(\varphi_2)$,
$\cos(\varphi_3)$ are distinct. Therefore, the matrix
$A=\cos^{\ell}(\varphi_j)_{\ell=0,1,2}^{j=1,2,3}$ is nonsingular.
Since $\cos(2\varphi)=2\cos^2(\varphi)-1$, the matrix
$\cos(\ell\varphi_j)_{\ell=0,1,2}^{j=1,2,3}$ is also nonsingular,
and this implies (4.1).

Next, by Lemma 4.1 (c), $\sin(\varphi_j)\neq0$ ($j=1,2,3$). 
Therefore, the matrix
$B=\sin(\varphi_j)\cos^{\ell}(\varphi_j)_{\ell=0,1,2}^{j=1,2,3}$
is nonsingular. Using the identities
$\sin(2\varphi)=2\sin(\varphi)\cos(\varphi)$,
$\sin(3\varphi)=2\sin(\varphi)(4\cos^2(\varphi)-1)$,
it follows that the matrix
$\sin(\ell\varphi_j)_{\ell=1,2,3}^{j=1,2,3}$ is also nonsingular.
This implies (4.2).
\enddemo

%%%%%%%%%%%%%%%%
%
% next, combine Lemmas 4,5,6 from old version
%
%%%%%%%%%%%%%%%%

\proclaim{Lemma 4.3} Let $z_1$ and $z_2$ be complex numbers.
We can associate with each $\chi\in C_q$
a non-negative real number $\lambda_\chi$ such that
$$
\aligned
z_1&=\sum_{\chi\in C_q} \lambda_\chi(\bc(a_2)-\bc(a_1)), \\
z_2&=\sum_{\chi\in C_q} \lambda_\chi(\bc(a_3)-\bc(a_2)).
\endaligned\tag{4.3}
$$
\endproclaim

\demo{Proof}  Write $z_j=u_j+iv_j$ ($j=1,2$), where $u_1,u_2,v_1,v_2$ are
real.
By Lemma 4.2, there is a character $\chi=\chi_0$ for which (4.1)
and (4.2) hold. Thus, we can find real numbers $\lambda_1$ and $\lambda_2$
such that
$$
\split
\lambda_1\Re (\chi_0(a_2)-\chi_0(a_1))+
\lambda_2\Re(\chi^2_0(a_2)-\chi^2_0(a_1))&=u_1/2, \\
\lambda_1\Re (\chi_0(a_3)-\chi_0(a_2))+
\lambda_2\Re(\chi^2_0(a_3)-\chi^2_0(a_2))&=u_2/2,
\endsplit
$$
and real numbers $\lambda_3$ and $\lambda_4$ such that
$$
\split
\lambda_3\Im (\chi^h_0(a_2)-\chi^h_0(a_1))+
\lambda_4\Im (\chi^k_0(a_2)-\chi^k_0(a_1))&=v_1/2, \\
\lambda_3\Im (\chi^h_0(a_3)-\chi^h_0(a_2))+
\lambda_4\Im (\chi^k_0(a_3)-\chi^k_0(a_2))&=v_2/2,
\endsplit
$$
By Lemma 4.1, the six characters $\chi_0,\chi_0^2,\chi_0^3$, $\bc_0,\bc_0^2,
\bc_0^3$ are distinct.
Now set $\mu_\chi=\lam_1$ for $\chi\in \{\chi_0, \bc_0\}$,
$\mu_\chi=\lam_2$ for $\chi\in \{ \chi_0^2, \bc_0^2 \}$, and $\mu_\chi=0$
for other characters.  Also, let
$\nu_{\chi_0^h}=\lam_3$, $\nu_{\bc_0^h}=-\lam_3$,
$\nu_{\chi_0^k}=\lam_4$, $\nu_{\bc_0^k}=-\lam_4$, and $\nu_\chi=0$ 
for other characters.  
Let $\theta_\chi = \mu_\chi + \nu_\chi$ for
each $\chi$.  Then (4.3) holds with $\lambda_\chi=\theta_\chi$ for each 
$\chi$, but it may occur that $\theta_\chi<0$ for some $\chi$.
However, by Lemma 4.1, $a_j \not\equiv 1 \pmod{q}$ for each $j$, so
$\sum_{\chi\in C_q} \chi(a_j) = -1$ for every $j$.  Thus, for any
real $y$, (4.3) holds
with $\lambda_\chi=\theta_\chi+y$ for each $\chi$.
\enddemo

\proclaim{Lemma 4.4} If $a_1,a_2,a_3$ do not satisfy the conditions of
Lemma 2.1, then for all $\tau>0$ and $\sigma>\frac12$,
there is a finite barrier for $D=(q,a_1,a_2,a_3)$, with each $B(\chi)$
consisting of numbers $\rho$ with $\Re \rho \le \sg$ and $\Im \rho>\tau$.
\endproclaim

\demo{Proof}
By Lemma 4.3, we can find such nonnegative $\nu^{(1)}_\chi$ and
$\nu^{(2)}_\chi$ that
$$
\aligned
 i&=\sum_\chi \nu^{(1)}_\chi(\overline\chi(a_2)-\overline\chi(a_1)),\\
-i&=\sum_\chi \nu^{(1)}_\chi(\overline\chi(a_3)-\overline\chi(a_2)),\\
 i&=\sum_\chi \nu^{(2)}_\chi(\overline\chi(a_2)-\overline\chi(a_1)),\\
 i&=\sum_\chi \nu^{(2)}_\chi(\overline\chi(a_3)-\overline\chi(a_2)).\\
\endaligned\tag{4.4}
$$
Fix small positive $\varepsilon>0$ and take a positive integer $Q$
and nonnegative integers $N^{(1)}_\chi$, $N^{(2)}_\chi$ for all characters
$\chi$ modulo $q$ such that
$|\nu^{(1)}_\chi-N^{(1)}_\chi/Q|<\varepsilon$,
$|\nu^{(2)}_\chi-N^{(2)}_\chi/Q|<\varepsilon$.
For some $\sigma_1 \in (\beta_1,\sigma]$ and large $\gamma>\tau$,
suppose that for all characters $\chi\in C_q$ and for $k=1,2$
the function $L(s,\chi)$
has a zero at $s=\sg_1+ ki\gamma$ of order $N^{(k)}_\chi$.
Suppose all other non-trivial zeros of $L$-functions
modulo $q$ have real part $\le \beta_1$.
Let $D_1(x)=\phi(q) (\pi_{q,a_1}(x)-\pi_{q,a_2}(x))$ and
$D_2(x)=\phi(q) (\pi_{q,a_2}(x)-\pi_{q,a_3}(x))$.
By Lemma 1.1 and (4.4), we have
$$
\frac{\log x}{x^{\sg_1}} D_1(x)=
\frac{Q}{2\gamma} (2\cos(\gamma\log x)+\cos(2\gamma\log x)+\varepsilon_1(x)
+ O(1/\gamma)),
$$
$$
\frac{\log x}{x^{\sg_1}} D_2(x)=
\frac{Q}{2\gamma} (-2\cos(\gamma\log x)+\cos(2\gamma\log x)+\varepsilon_2(x)
+ O(1/\gamma)),
$$
where the functions $\varepsilon_1(x)$, $\varepsilon_2(x)$ are uniformly
small if $\varepsilon$ is small. Taking into account that
$\min(2\cos u+\cos2u,-2\cos u+\cos2u)\le-1$ for all $u$, we obtain that
for large $x$, $\pi_{q,a_1}(x) > \pi_{q,a_2}(x) > \pi_{q,a_3}(x)$ does
not occur.
\qed\enddemo

%%%%%%%%%%%%%%%%%%%%%%%%%%%%%%%%%%%%%%%%%%%%%%%%%%%%%%%%%%%%%%%%%%%%%%%%
%
%
\head 5. A barrier satisfying GSH$_q$ \endhead
%
%
%%%%%%%%%%%%%%%%%%%%%%%%%%%%%%%%%%%%%%%%%%%%%%%%%%%%%%%%%%%%%%%%%%%%%%%%

The construction of this barrier is modeled on the construction in \S 2.
For one character, $B(\chi)$ is infinite, the number of elements of $B(\chi)$
with imaginary part $\le T$ growing like $\sqrt{T}$.  By altering the
parameters in the construction, we can create barriers with $\sqrt{T}$
replaced by $T^{\epsilon}$ for any fixed $\epsilon$.
Assume that for some relabeling of $a_1,a_2,a_3$,
 there are two characters $\chi_1$,$\chi_2$ satisfying
$$
\chi_1(a_1)=\chi_1(a_2) \ne \chi_1(a_3), \qquad \chi_2(a_1) \ne 
\chi_2(a_2). \tag{5.1}
$$
Suppose that $\frac12 \le \beta < \sigma_2 <\sigma_1$, that $t$ is
large and that $L(s,\chi_1)$ has a simple zero at $s=\sigma_1+it$.  Suppose
that $L(s,\chi_2)$ has simple zeros at the points $s=\rho_j$ ($j=1,2,\ldots$),
where $\rho_j=\sigma_2-\delta_j+i
\g_j$, $\del_j>0$, $\g_j>0$, $\del_j \to 0$ and $\g_j\to \infty$ as
$j\to\infty$, and
$\sum 1/\g_j < \infty$. Also, suppose the numbers $t,\g_1,\g_2,\ldots$
are linearly independent over $\Bbb Q$.
 Define
$$
Z=\bc_1(a_2)-\bc_1(a_3), \qquad W=\bc_2(a_2)-\bc_2(a_1).
$$
By (5.1), $Z\ne 0$ and $W\ne 0$.  Also define
$$
\a = -\frac{1}{\pi} \( \tan^{-1} \frac{\sigma_1}{t} + \arg Z \), \quad
\beta = \frac{\arg W}{2\pi}-\frac14.
$$
Let $H$ be the set of integers $h$ such that $\| h\a + \b \| \le \frac15$.
Since the number of  possibilities for $Z$ is finite, if $t$ is large then
$$
\frac{1}{10t} \le \| \a \| \le \frac12 -\frac{1}{10t}.
$$
It follows that in every set of $\lfloor 10t \rfloor+1$ consecutive integers,
one of them is in $H$.  As in section 2, define
$$
D_1(x):=\phi(q)(\pi_{q,a_1}(x) - \pi_{q,a_2}(x)), \qquad
D_2(x):=\phi(q)(\pi_{q,a_3}(x) - \pi_{q,a_2}(x)).
$$
Suppose $x$ is sufficiently large, and for brevity write $u=\log x$.
By Lemma 1.1 and our hypotheses,
$$
D_2(x) = \frac{2 x^{\sigma_1}}{u} \left[
\Re \(\frac{e^{itu}}{\sigma_1 + it} Z\) + O\pfrac{1}{u} \right] \tag{5.2}
$$
and
$$
\aligned
D_1(x) &= \frac{2 x^{\sigma_2}}{u} \sum_{\g_j \le x} \left[ \Re
  \( \frac{e^{(-\del_j+i\g_j)u}}{\sg_2-\del_j+i\g_j}W \)
   + O\pfrac{e^{-\del_j u}}
  {\g_j^2 u} \right] + O(x^\beta\log^2 x) \\
&= \frac{2 x^{\sigma_2}}{u} \sum_{\g_j \le x} \left[ \Re B_j + O\pfrac{
  e^{-\del_j u}}{\g_j^2} \right] + O(x^\beta\log^2 x),
\endaligned\tag{5.3}
$$
where
$$
B_j = W \frac{e^{(-\del_j+i\g_j)u}}{i\g_j}. \tag{5.4}
$$
By assumption, $\sum_j |B_j| \ll 1$, thus $D_1(x) \ll x^{\sg_2}/u$.
Modulo $2\pi$,
$$
\arg \frac{e^{itu}}{\sigma_1 + it} Z \equiv tu-\tan^{-1} \frac{t}{\sigma_1}
+\arg Z \equiv tu-\frac{\pi}2-\pi \a.
$$
By (5.2), when $\| tu/\pi - \a \| \ge u^{-0.9}$, $D_2(x) \gg x^{\sg_1}/(\log x)
^{1.9}$, and thus for these $x$ either $\pi_{q,a_3}(x)$ is the largest or
smallest of the three functions.  Next assume that
$$
\| tu/\pi - \a \| \le u^{-0.9}.
$$
We choose $\del_j$ and $\g_j$ as follows: $0 < \del_j < \sg_2-\beta$,
$j^{-3} \ll \del_j \ll j^{-3}$, $\g_j=2th_j+O(j^{-10})$, where
for $j\ge 10t$ we have $h_j\in H$, $h_{j+1}>h_j$ and $j^2 \le h_j \le j^2+j$.
With these choices,
$$
\sum_{j=1}^\infty \frac{e^{-\del_j u}}{\g_j^2} \ll e^{-u^{1/4}} \sum_{j\le
u^{1/4}} 1/j^4 + \sum_{j>u^{1/4}}1/j^4 \ll u^{-3/4}
$$
and
$$
\sum_{j<u^{1/4}\text{ or }j>u^{2/5}}  \frac{e^{-\del_j u}}{\g_j} \ll 
e^{-u^{1/4}}+u^{-2/5} \ll u^{-2/5}.
$$
Thus, by (5.3) and (5.4),
$$
D_1(x) = \frac{2 x^{\sigma_2}}{u} \left[ 
\sum_{u^{1/4} \le j \le u^{2/5}} \Re B_j + O(u^{-0.4}) \right]. \tag{5.5}
$$
Suppose $u^{1/4} \le j\le u^{2/5}$.  Since $h_j \in H$, we have
\def\lnm{\left\|}
\def\rnm{\right\|}
$$
\split
\lnm \frac{1}{2\pi} \arg B_j \rnm &= \lnm \frac{1}{2\pi}\(\arg W+\g_j u -
\frac{\pi}{2}\) \rnm \\
&=\lnm \beta + \frac{ut}{\pi} h_j + O(u^{-3/2}) \rnm \\
&= \lnm \beta +h_j \a + O(u^{-0.1})\rnm \le 0.21
\endsplit
$$
for large $u$.  Hence $\Re B_j \ge |B_j| \cos(0.42\pi) \ge \frac15|B_j|$.
Therefore,
$$
\sum_{u^{1/4} \le j\le u^{2/5}} \Re B_j \gg \sum_{u^{1/3} \le j \le 2u^{1/3}}
\frac{1}{\g_j} \gg u^{-1/3}.
$$
It follows from (5.5)
that for $u$ large and $\| \frac{ut}{\pi}-\a \| \le u^{-0.9}$ that
$$
D_1(x) \ge \frac{c x^{\sg_2}}{(\log x)^{4/3}}
$$
where $c>0$ depends on $q$, $t$ and $W$.  This implies
that the inequality $\pi_{q,a_2}(x) > \pi_{q,a_3}(x) > \pi_{q,a_1}(x)$
does not occur for large $x$.

\enddocument

\Refs
\refstyle{A}
\widestnumber\key{BFHR}

\ref\key{BH} \by C. Bays and R. H.\ Hudson \paper Details of the first
region of integers $x$
with $\pi_{3,2}(x)< \pi_{3,1}(x)$ \jour Math.\ Comp.\ \vol 32 \yr 1978
\pages 571--576 \endref

\ref\key{Ch}  \by P.\ L.\ Chebyshev \paper  Lettre de M. le professeur
Tch\'ebychev \`a M. Fuss, sur un
nouveau th\'eoreme r\'elatif aux nombres premiers contenus dans la formes
$4n+1$ et
$4n+3$ \jour  Bull.\ de la Classe phys.-math. de l\rq Acad.\ Imp.\ des
Sciences St.\ Petersburg \vol 11
\yr 1853 \page 208 \endref

\ref\key{D} \by H. Davenport \book Multiplicative Number Theory,
3rd ed., Graduate Texts in Mathematics vol. 74 \publ Springer-Verlag
\publaddr New York-Berlin \yr 2000
\endref

\ref\key{K1} \by J. Kaczorowski \paper A contribution to the Shanks-R\'enyi
 race problem \jour  Quart. J. Math., Oxford Ser. (2) \vol 449 \yr 1993
\pages 451--458 \endref

\ref\key{K2} \bysame \paper On the Shanks-R\'enyi race mod $5$
\jour J. Number Theory \vol 50 \yr 1995 \pages 106--118
\endref

\ref\key{K3} \bysame \paper On the Shanks-R\'enyi race problem
\jour Acta Arith. \vol 74 \yr 1996 \pages 31--46
\endref

\ref\key{KT1} \by S. Knapowski and P. Tur\'an \paper
Comparative prime number theory I.
\jour Acta. Math. Sci. Hungar. \vol 13 \yr 1962 \pages 299-314
\moreref II.\vol 13 \yr 1962 \pages 315--342 
\moreref III.\vol 13 \yr 1962 \pages 343--364
\moreref IV. \vol 14 \yr 1963 \pages 31--42 
\moreref V. \vol 14 \yr 1963\pages 43--63
\moreref VI.\vol 14 \yr 1963 \pages 65--78 
\moreref VII.\vol 14 \yr 1963 \pages 241--250 
\moreref VIII. \vol 14 \yr 1963 \pages 251--268 \endref

\ref\key{KT2} \by S. Knapowski and P. Tur\'an \paper
 Further developments in the
comparative prime-number theory. I. \jour Acta Arith. \vol 9 \yr 1964
\pages 23--40 
\moreref II. \vol 10 \yr 1964 \pages 293--313
\moreref III. \vol 11 \yr 1965 \pages 115--127
\moreref IV. \vol 11 \yr 1965 \pages 147--161
\moreref V. \vol 11 \yr 1965 \pages 193--202
\moreref VI. \vol 12 \yr 1966 \pages 85--96
\endref

%\ref\key{L} \by R. S. Lehman \paper On the difference $\pi(x)-\li(x)$
%\jour Acta Arith. \vol 11 \yr 1966 \pages 397--410 \endref
%
\ref\key{L} \by J.E. Littlewood \paper
Sur la distribution des nombres premiers
\jour C. R. Acad. des sciences Paris
\vol 158 \yr 1914 \pages 1869--1872 \endref

\ref\key{RS}\by  M. Rubinstein and P. Sarnak \paper Chebyshev's Bias
\jour J. Exper. Math. \vol 3 \yr 1994 \pages 173--197 \endref

\ref\key{Ru} \by R. Rumely \paper Numerical computations concerning the ERH
\jour Math. Comp. \vol 61 \yr 1993 \pages 415--440 \endref

\ref\key{Sh} \by D. Shanks \paper Quadratic residues and the distribution of
primes \jour Math. Comp. \vol 13 \yr 1959 \pages  272--284 \endref

\endRefs
\bye